\documentclass[12pt]{article}
\usepackage{amssymb,amsmath,latexsym}
\allowdisplaybreaks

\usepackage{color}

\begin{document}

\begin{doublespace}

\def\1{{\bf 1}}
\def\ind{{\bf 1}}
\def\nn{\nonumber}
\newcommand{\I}{\mathbf{1}}
\def\X {{\mathfrak X}}
\def\sA {{\cal A}} \def\sB {{\cal B}} \def\sC {{\cal C}}
\def\sD {{\cal D}} \def\sE {{\cal E}} \def\sF {{\cal F}}
\def\sG {{\cal G}} \def\sH {{\cal H}} \def\sI {{\cal I}}
\def\sJ {{\cal J}} \def\sK {{\cal K}} \def\sL {{\cal L}}
\def\sM {{\cal M}} \def\sN {{\cal N}} \def\sO {{\cal O}}
\def\sP {{\cal P}} \def\sQ {{\cal Q}} \def\sR {{\cal R}}
\def\sS {{\cal S}} \def\sT {{\cal T}} \def\sU {{\cal U}}
\def\sV {{\cal V}} \def\sW {{\cal W}} \def\sX {{\cal X}}
\def\sY {{\cal Y}} \def\sZ {{\cal Z}}

\def\bA {{\mathbb A}} \def\bB {{\mathbb B}} \def\bC {{\mathbb C}}
\def\bD {{\mathbb D}} \def\bE {{\mathbb E}} \def\bF {{\mathbb F}}
\def\bG {{\mathbb G}} \def\bH {{\mathbb H}} \def\bI {{\mathbb I}}
\def\bJ {{\mathbb J}} \def\bK {{\mathbb K}} \def\bL {{\mathbb L}}
\def\bM {{\mathbb M}} \def\bN {{\mathbb N}} \def\bO {{\mathbb O}}
\def\bP {{\mathbb P}} \def\bQ {{\mathbb Q}} \def\bR {{\mathbb R}}
\def\bS {{\mathbb S}} \def\bT {{\mathbb T}} \def\bU {{\mathbb U}}
\def\bV {{\mathbb V}} \def\bW {{\mathbb W}} \def\bX {{\mathbb X}}
\def\bY {{\mathbb Y}} \def\bZ {{\mathbb Z}}
\def\R {{\mathbb R}} \def\RR {{\mathbb R}} \def\H {{\mathbb H}}
\def\n{{\bf n}} \def\Z {{\mathbb Z}}

\newcommand{\expr}[1]{\left( #1 \right)}
\newcommand{\cl}[1]{\overline{#1}}
\newtheorem{thm}{Theorem}[section]
\newtheorem{lemma}[thm]{Lemma}
\newtheorem{defn}[thm]{Definition}
\newtheorem{prop}[thm]{Proposition}
\newtheorem{corollary}[thm]{Corollary}
\newtheorem{remark}[thm]{Remark}
\newtheorem{example}[thm]{Example}
\numberwithin{equation}{section}
\def\ee{\varepsilon}
\def\qed{{\hfill $\Box$ \bigskip}}
\def\NN{{\mathcal N}}
\def\AA{{\mathcal A}}
\def\MM{{\mathcal M}}
\def\BB{{\mathcal B}}
\def\CC{{\mathcal C}}
\def\LL{{\mathcal L}}
\def\DD{{\mathcal D}}
\def\FF{{\mathcal F}}
\def\EE{{\mathcal E}}
\def\QQ{{\mathcal Q}}
\def\SS{{\mathcal S}}
\def\RR{{\mathbb R}}
\def\R{{\mathbb R}}
\def\L{{\bf L}}
\def\K{{\bf K}}
\def\S{{\bf S}}
\def\A{{\bf A}}
\def\E{{\mathbb E}}
\def\F{{\bf F}}
\def\P{{\mathbb P}}
\def\N{{\mathbb N}}
\def\eps{\varepsilon}
\def\wh{\widehat}
\def\wt{\widetilde}
\def\pf{\noindent{\bf Proof.} }
\def\pff{\noindent{\bf Proof} }
\def\cp{\mathrm{Cap}}

\title{\Large \bf Martin boundary of unbounded sets  for purely discontinuous Feller processes}

\author{{\bf Panki Kim}\thanks{This work was  supported by the National Research Foundation of
Korea(NRF) grant funded by the Korea government(MSIP) (No. NRF-2015R1A4A1041675)
}
\quad {\bf Renming Song\thanks{Research supported in part by a grant from
the Simons Foundation (208236)}} \quad and
\quad {\bf Zoran Vondra\v{c}ek}
\thanks{Research supported in part by the Croatian Science Foundation under the project 3526}
}

\date{}

\maketitle

\begin{abstract}
In this paper, we study the Martin kernels of general open sets associated with inaccessible
points for a large class of purely discontinuous Feller processes in metric measure spaces.
\end{abstract}



\section{Introduction and setup}\label{s:intro}

This paper is a companion of \cite{KSVp2} and here we continue our study of the Martin boundary
of Greenian open sets with respect to purely discontinuous Feller processes in metric measure spaces.
In \cite{KSVp2}, we have shown that
(1) if $D$ is a Greenian open set and $z_0\in \partial D$ is accessible from $D$, then the Martin kernel of
$D$ associated with $z_0$ is a minimal harmonic function; (2) if $D$ is an unbounded Greenian open set and
$\infty$ is accessible from $D$, then the Martin kernel of $D$ associated with $\infty$ is a minimal harmonic function.
The goal of this paper is to study the Martin kernels of $D$ associated with inaccessible boundary points
of $D$, including $\infty$.

The background and recent progress on the Martin boundary is explained in the companion paper
\cite{KSVp2}.
Martin kernels of bounded open sets $D$ associated with both accessible and inaccessible
boundary points of $D$ have been studied in the recent preprint \cite{JK}. In this paper,
we are mainly concerned with the Martin kernels of unbounded open sets associated with $\infty$
when $\infty$ is inaccessible from $D$. For completeness, we also spell out
some of the details of the argument for dealing with the Martin kernels of unbounded open sets associated with
inaccessible boundary points of $D$.
To accomplish our task of studying the Martin kernels of general open sets,
we follow the ideas of \cite{BKK, KL} and first study the oscillation reduction of ratios of positive harmonic functions.
In the case of isotropic $\alpha$-stable processes, the oscillation reduction at infinity and Martin kernel
associated with $\infty$ follow easily from the corresponding results at finite boundary points
by using the sphere inversion and Kelvin transform.
For the general processes dealt with in this paper, the Kelvin transform method does not apply.

Now we describe the setup of this paper which is the same as that of \cite{KSVp2}
and then give the main results of this paper.

Let $(\X, d, m)$ be a metric measure space with a countable base
 such that all bounded closed sets are compact and the measure $m$ has full support. For $x\in \X$ and $r>0$, let $B(x,r)$ denote the ball centered at $x$ with radius $r$. Let $R_0\in (0,\infty]$ be the localization radius such that $\X\setminus B(x,2r)\neq \emptyset$ for all $x\in \X$ and all $r<R_0$.

Let $X=(X_t, \sF_t, \P_x)$ be a Hunt process on $\X$. We will assume the following

\smallskip
\noindent
\textbf{Assumption A:}
$X$ is a Hunt process admitting a strong dual process $\widehat{X}$ with respect to the measure $m$ and $\widehat{X}$ is also a Hunt process.
The transition semigroups $(P_t)$ and $(\widehat{P}_t)$ of $X$ and $\widehat{X}$ are both Feller and strong Feller.  Every semi-polar set of $X$ is polar.

\smallskip
In the sequel, all objects related to the dual process $\widehat{X}$ will be denoted by a hat.
We first recall that a set is polar (semi-polar, respectively) for $X$ if and only if it is polar (semi-polar, respectively) for $\wh X$.

 If $D$ is an open subset of $\X$ and $\tau_D=\inf\{t>0:\, X_t\notin D\}$ the exit time from $D$, the killed process $X^D$ is defined by $X_t^D=X_t$ if $t<\tau_D$ and $X_t^D=\partial$
 where $\partial$ is an extra point added to $\X$.
 Then, under assumption {\bf A},
 $X^D$ admits a unique (possibly infinite) Green function (potential kernel) $G_D(x,y)$ such that for every non-negative Borel function $f$,
$$
G_D f(x):=\E_x \int_0^{\tau_D}f(X_t)dt=\int_D G_D(x.y)\, m(dy)\,  ,
$$
and $G_D(x,y)=\widehat{G}_D(y,x)$, $x,y\in D$, with $\widehat{G}_D(y,x)$ the Green function of $\widehat{X}^D$.
It is assumed throughout the paper that $G_D(x,y)=0$ for $(x,y)\in (D\times D)^c$.
We also note that the killed process $X^D$ is strongly Feller, see e.g.~the first part of the proof of Theorem on \cite[pp.~68--69]{Chu}.

Let $\partial D$ denote the boundary of the open set $D$ in the topology of $\X$.
Recall that $z\in \partial D$ is said to be regular for $X$ if $\P_z(\tau_D =0)=1$ and irregular otherwise.
We will denote the set of regular points of $\partial D$ for $X$ by $D^{\mathrm{reg}}$ (and the set of regular points of $\partial D$ for $\wh X$ by $\wh D^{\mathrm{reg}}$).
It is well known that the set of irregular points is semipolar, hence polar under {\bf A}.

 Suppose that $D$ is Greenian, that is, the Green function $G_D(x,y)$ is finite away from the diagonal.
Under this assumption, the killed process $X^D$ is transient (and strongly Feller). In particular, for every bounded Borel function $f$ on $D$, $G_D f$ is continuous.

The process $X$, being a Hunt process, admits a L\'evy system $(J,H)$ where $J(x,dy)$
is a kernel on $\X$ (called the L\'evy kernel of $\X$),
and $H=(H_t)_{t\ge 0}$ is a positive continuous additive functional of $X$. We assume that $H_t=t$
so that for every function $f:\X\times \X \to [0,\infty)$ vanishing on the diagonal
and every stopping time $T$,
$$
\E_x \sum_{0<s\le T} f(X_{s-}, X_s)=\E_x \int_0^T f(X_s,y)J(X_s,dy) ds\, .
$$
Let $D\subset \X$ be a Greenian  open set. By replacing $T$ with $\tau_D$ in the displayed formula above and taking $f(x,y)=\1_D(x)\1_A(y)$ with $A\subset \overline{D}^c$, we get that
\begin{equation}\label{e:exit-distribution}
\P_x(X_{\tau_D}\in A, \tau_D <\zeta)=\E_x \int_0^{\tau_D} J(X_s, A) ds= \int_D G_D(x,y)J(y,A)m(dy)\,
\end{equation}
where $\zeta$ is the life time of $X$.
Similar formulae hold for $\wh{X}$ and $\wh{J}(x,dy)m(dx)=J(y,dx)m(dy)$.
\medskip
\noindent
\textbf{Assumption C:} The L\'evy kernels of $X$ and $\wh{X}$ have the form $J(x,dy)=j(x,y)m(dy)$, $\wh{J}(x,dy)=\widehat{j}(x,y)m(dy)$, where $j(x,y)=\widehat{j}(y,x)>0$ for all $x,y\in \X$, $x\neq y$.

\medskip

We will always assume that Assumptions \textbf{A} and \textbf{C} hold true.

\medskip

In the next assumption, $z_0$ is a point in $\X$ and $R\le R_0$.

\medskip

\noindent
\textbf{Assumption C1}$(z_0, R)$:
For all $0<r_1<r_2<R$, there exists a constant $c=c(z_0, r_2/r_1)>0$ such that for all $x\in B(z_0,r_1)$ and all $y\in \X\setminus B(z_0, r_2)$,
$$
 c^{-1} j(z_0,y)\le j(x,y)\le c j(z_0,y), \qquad  c^{-1} \widehat{j}(z_0,y)\le \widehat{j}(x,y)\le c \widehat{j}(z_0,y).
$$

\medskip

In the next assumption we require that the localization radius $R_0=\infty$ and that $D$ is unbounded.
Again, $z_0$ is a point in $\X$.

\medskip

\noindent
\textbf{Assumption C2}$(z_0, R)$:
For all $R\le r_1<r_2< \infty$, there exists a constant $c=c(z_0, r_2/r_1)>0$ such that for all $x\in B(z_0,r_1)$ and all $y\in \X\setminus B(z_0, r_2)$,
$$
 c^{-1} j(z_0,y)\le j(x,y)\le c j(z_0,y), \qquad  c^{-1} \widehat{j}(z_0,y)\le \widehat{j}(x,y)\le c \widehat{j}(z_0,y).
$$

\medskip
We \emph{define}  the Poisson kernel of $X$ on  an open set $D\in \X$  by
\begin{equation*}
P_D(x,z)=\int_D G_D(x,y)j(y,z) m(dy), \qquad x\in D, z\in D^c .
\end{equation*}
By \eqref{e:exit-distribution}, we see that $P_D(x,\cdot)$ is the density of the exit distribution of $X$ from $D$
restricted to $\overline{D}^c$:
$$
\P_x(X_{\tau_D}\in A, \tau_D <\zeta)=\int_A P_D(x,z) m(dz), \qquad A\subset \overline{D}^c .
$$

Recall that $f:\X\to [0,\infty)$ is regular harmonic in $D$ with respect to $X$ if
$$
f(x)=\E_x[f(X_{\tau_D}), \tau_D<\zeta]\, , \quad \textrm{for all }x\in D\,  ,
$$
and it is harmonic in $D$ with respect to $X$ if for every relatively compact open $U\subset \overline{U}\subset D$,
$$
f(x)=\E_x[f(X_{\tau_U}), \tau_D<\zeta]\, , \quad \textrm{for all }x\in U\,  .
$$
Recall also that $f:D\to [0,\infty)$ is harmonic in $D$ with respect to $X^D$ if for every relatively compact open $U\subset \overline{U}\subset D$,
$$
f(x)=\E_x[f(X^D_{\tau_U}), \tau_U < \zeta]\, , \quad \textrm{for all }x\in U\,  .
$$

The next pair of assumptions is about an approximate factorization of
positive harmonic functions. This approximate factorization plays a crucial role in
proving the oscillation reduction.
The first one is an approximate factorization of harmonic functions at a finite boundary point.

\noindent
\textbf{Assumption F1}$(z_0, R)$:
Let $z_0\in \X$ and $R\le R_0$. For any $\frac{1}{2} < a < 1$, there exists $C(a)=C(z_0, R, a)\ge 1$ such that for
every $r\in (0, R)$, every open set $D\subset B(z_0,r)$, every nonnegative function $f$ on $\X$
which is  regular harmonic in $D$ with respect to $X$ and vanishes in
$B(z_0, r) \cap ( \overline{D}^c \cup D^{\mathrm{reg}})$,
and all $x\in D\cap B(z_0,r/8)$\,,
\begin{eqnarray}\label{e:af-1}
\lefteqn{C(a)^{-1}\E_x[\tau_{D}] \int_{\overline{B}(z_0,ar/2)^c}  j(z_0,y) f(y)m(dy)  }\nonumber \\
&&\le f(x) \le C(a)\E_x[\tau_{D}]\int_{\overline{B}(z_0,ar/2)^c}   j(z_0,y) f(y)m(dy).
\end{eqnarray}

In the second assumption we require that the localization radius $R_0=\infty$ and that $D$ is unbounded.

\noindent
\textbf{Assumption F2}$(z_0, R)$:
Let $z_0\in \X$ and $R>0$. For any $1 < a  < 2$, there exists  $C(a)=C(z_0, R, a)\ge 1$ such that for
every $r\ge R$, every open set $D\subset \overline{B}(z_0,r)^c$,
every nonnegative function $f$ on $\X$ which is  regular harmonic in $D$
with respect to $X$ and vanishes
on
$\overline{B}(z_0, r)^c  \cap ( \overline{D}^c \cup  D^{\mathrm{reg}})$,
and all $x\in D\cap \overline{B}(z_0,8r)^c$,
\begin{eqnarray}\label{e:af-2}
\lefteqn{C(a)^{-1}\,  P_{D}   (x,z_0) \int_{B(z_0, 2ar)}  f(z)m(dz)}\nonumber\\
&&\le f(x) \le
C(a)\,   P_{D}   (x,z_0) \int_{B(z_0, 2ar)}  f(z)m(dz).
\end{eqnarray}

Let $D\subset \X$ be an open set. A point $z\in \partial D$ is said to be accessible from $D$ with respect to $X$ if
\begin{equation}\label{e:accessible-finite}
P_D(x,z)=\int_D G_D(x,w)j(w,z) m(dw) = \infty \quad \text{ for all } x \in D\, ,
\end{equation}
and inaccessible otherwise.

In case $D$ is unbounded we say that $\infty$ is accessible from $D$ with respect to $X$ if
\begin{equation}\label{e:accessible-finite2}
\E_x \tau_D =\int_D G_D(x,w) m(dw)=  \infty \quad \text{ for all } x \in D
\end{equation}
 and inaccessible otherwise.
 The notion of accessible and inaccessible points was introduced in \cite{BKuK}.

In \cite{KSVp2}, we have discussed the oscillation reduction and Martin boundary points at accessible points, and
showed that the Martin kernel associated with an accessible point is a minimal harmonic function.
As in \cite{KSVp2},
the main tool in studying the Martin kernel associated with inaccessible
points is the oscillation reduction at inaccessible points. To prove the oscillation reduction
at inaccessible points,
we need to assume one of the following additional conditions
on the asymptotic behavior of the L\'evy kernel:

\noindent
\textbf{Assumption E1}$(z_0, R)$: For every $r \in (0, R)$,
$$
\lim_{d(z_0, y)  \to 0}\sup_{z: d(z_0, z)>r}\frac{j(z, z_0)}{j( z, y)}=\lim_{d(y, z_0)  \to 0}\inf_{z:d(z_0, z)>r}\frac{j(z, z_0)}{j(z, y)}=1.
$$

\noindent
\textbf{Assumption E2}$(z_0, R)$: For every $r>R$,
$$
\lim_{d(z_0, z) \to \infty} \sup_{y: d(z_0, y) < r}
\frac{j(z, z_0)}{j(z, y)}=\lim_{d(z_0, z) \to \infty} \inf_{y: d(z_0, y) < r}
\frac{j(z, z_0)}{j(z, y)}
=1.
$$

Combining Theorems \ref{t:oscillation-reduction-yI} and \ref{t:reduction} below for inaccessible points with the results in  \cite{KSVp2} for accessible ones, we have the following,
which is the first main result of this paper.

\begin{thm}\label{t:main-mb0}
Let $D\subset \X$ be an open set.
(a) Suppose that $z_0\in \partial D$.
Assume that there exists $R\le R_0$ such that
{\bf C1}$(z_0, R)$  and {\bf E1}$(z_0, R)$ hold,
and that $\widehat{X}$ satisfies \textbf{F1}$(z_0, R)$.
Let $r  \le R$ and let $f_1$ and $f_2$ be nonnegative functions on $\X$ which are regular harmonic in $D\cap {B}(z_0, r)$ with respect to $\widehat{X}$ and vanish  on
$B(z_0, r) \cap ( \overline{D}^c \cup \wh D^{\mathrm{reg}})$.
 Then the limit
$$
\lim_{D\ni x\to z_0}\frac{f_1(x)}{f_2(x)}
$$
exists and is finite.

\noindent
(b)
Suppose that $R_0=\infty$ and $D$ is an unbounded subset of $\X$.
Assume that there is a point $z_0\in\X$ such that
{\bf C2}$(z_0, R)$ and {\bf E2}$(z_0, R)$ hold,
and that  $\widehat{X}$ satisfies \textbf{F2}$(z_0, R)$ for some $R>0$.
Let $r > R$ and let $f_1$ and $f_2$ be nonnegative functions on $\X$ which are regular harmonic in $D\cap \overline{B}(z_0, r)^c$ with respect to $\widehat{X}$ and vanish  on
$\overline{B}(z_0, r)^c \cap ( \overline{D}^c \cup \wh D^{\mathrm{reg}}) $.
Then the limit
$$
\lim_{D\ni x\to \infty}\frac{f_1(x)}{f_2(x)}
$$
exists and is finite.
\end{thm}

For $D\subset \X$, let $\partial_M D$ denote the Martin boundary of $D$ with respect to $X^D$ in the sense of Kunita-Watanabe \cite{KW}, see Section 3 for more details. 
A point $w\in \partial_M D$ is said to be minimal if the Martin kernel $M_D(\cdot, w)$ is a minimal
harmonic function with respect to $X^D$. We will use $\partial_m D$ to denote the
minimal Martin boundary of $D$ with respect to $X^D$.
A point $w\in \partial_M D$ is said to be a \emph{finite Martin boundary point}
if there exists a bounded (with respect to the metric $d$) sequence $(y_n)_{n\ge 1}\subset D$
converging to $w$ in the Martin topology.
A point $w\in \partial_M D$ is said to be an \emph{infinite Martin boundary point}
if there exists an unbounded (with respect to the metric $d$) sequence $(y_n)_{n\ge 1}\subset D$
converging to $w$ in the Martin topology. We note that these two definitions do not rule out the possibility that a point $w\in \partial_M D$ is at the same time finite and infinite Martin boundary point. We will show in Corollary \ref{c:finite-not-infinite}(a) that under appropriate and natural assumptions this cannot happen. 
A point $w\in \partial_MD$ is said to be associated
with $z_0\in \partial D$ if there is a sequence $(y_n)_{n\ge 1}\subset D$ converging to $w$
in the Martin topology and to $z_0$ in the topology of $\X$. The set of Martin
boundary points associated with $z_0$ is denoted by $\partial_M^{z_0} D$.
A point $w\in \partial_MD$ is said to be associated with $\infty$ if $w$ is an infinite Martin boundary point.
The set of Martin boundary points associated with $\infty$ is denoted by $\partial_M^{\infty} D$.
$\partial^f_M D$ and $\partial^f_m D$ will be used to denote the finite part of the Martin
boundary and minimal boundary respectively.
Note that $\partial_M^{\infty} D$ is the set of infinite Martin boundary points.

 Recall that we denote the set of regular points of $\partial D$ for $X$ by $D^{\mathrm{reg}}$.
 Here is our final assumption.

  \medskip

 \noindent
\textbf{Assumption G}:
$\lim_{D\ni x\to z}G_D(x,y)=0$ for every $z\in D^{\mathrm{reg}}$ and every $y\in D$.

 \medskip

  \noindent

From Theorem \ref{t:main-mb0} and the results in \cite{KSVp2}, we have the following.
\begin{thm}\label{t:main-mb3}
Let $D\subset \X$ be an open set.
(a) Suppose that $z_0\in \partial D$.
Assume that there exists $R\le R_0$ such that {\bf C1}$(z_0, R)$
and {\bf E1}$(z_0, R)$ hold, and that
$\widehat{X}$ satisfies {\bf F1}$(z_0, R)$.
Then there is only one Martin boundary point associated with $z_0$.

\noindent
(b) Assume further that
{\bf G} holds,
${X}$ satisfies {\bf F1}$(z_0, R)$, and that for all $r\in (0,R]$,
\begin{align}\label{e:nGassup1}
 \sup_{x\in D\cap B(z_0,r/2)} \sup_{y \in \X \setminus B(z_0, r)} \max(G_D(x, y), \widehat{G}_D(x, y))=:c(r) <\infty,
 \end{align}
and in case of unbounded $D$, for  $r \in (0, r_0]$,
$$
\lim_{x \to \infty} G_D(x, y)=0 \quad \text{for all } y \in D \cap B(z_0, r)\, .
$$
Then the Martin boundary point associated with $z_0\in \partial D$ is minimal if and only if $z_0$ is accessible from $D$ with respect to $X$.
\end{thm}

\begin{corollary}\label{c:main-mb3}
Suppose that the assumptions of Theorem \ref{t:main-mb3}(b) are satisfied for all $z_0\in \partial D$ (with $c(r)$ in \eqref{e:nGassup1} independent of $z_0$). Suppose
further that, for any inaccessible point $z_0\in \partial D$, $\lim_{D \ni x\to z_0}j(x, z_0)=\infty$.

\noindent
(a)  Then the finite part of the Martin boundary $\partial_M D$ can be identified
with $\partial D$.

\noindent
(b) If $D$ is bounded, then $\partial D$ and $\partial_M D$ are homeomorphic.
\end{corollary}

\begin{thm}\label{t:main-mb4}
(a) Suppose that $R_0=\infty$ and $D$ is an unbounded open subset of $\X$.
If there is a point $z_0\in \X$ such that  {\bf C2}$(z_0, R)$
and {\bf E2}$(z_0, R)$ hold,
and $\wh{X}$ satisfies {\bf F2}$(z_0, R)$,
then there is only one Martin boundary point associated with $\infty$.

\noindent
(b)  Assume further that
{\bf G} holds, ${X}$ satisfies {\bf F2}$(z_0, R)$, and that
for all $r\ge R$,
\begin{align}\label{e:nGassup1-infty}
 \sup_{x\in D\cap B(z_0,r/2)} \sup_{y \in \X \setminus B(z_0, r)} \max(G_D(x, y), \widehat{G}_D(x, y))=:c(r) <\infty
 \end{align}
 and
\begin{align}
\label{e:G111}
\lim_{x \to \infty} G_D(x, y)=0 \quad \text{for all } y \in D.
\end{align}
 Then the Martin boundary point
associated with $\infty$ is minimal if and only if $\infty$ is accessible from $D$.
\end{thm}

\begin{corollary}\label{c:finite-not-infinite} Let $R_0=\infty$ and $D\subset \X$ be unbounded. Suppose that the assumptions of Theorem \ref{t:main-mb3}(b) are satisfied for all $z_0\in \partial D$ (with $c(r)$ in \eqref{e:nGassup1} independent of $z_0$) and that the assumptions of 
Theorem \ref{t:main-mb4}(a) and (b) are satisfied. Then 

\noindent
(a) $\partial_M^f D\cap \partial_M^{\infty}D=\emptyset$.

\noindent 
(b) Suppose that, for any inaccessible point $z_0\in \partial D$, $\lim_{D \ni x\to z_0}j(x, z_0)=\infty$. Then the Martin boundary $\partial_M D$ is homeomorphic with the one-point compactification of $\partial D$.
\end{corollary}

 In case when $X$ is an isotropic stable process, Theorems \ref{t:main-mb3} and \ref{t:main-mb4} were proved in \cite{BKK}.

In Section \ref{s:inaccessible} we provide the proof of Theorem \ref{t:main-mb0} for inaccessible points. Section \ref{s:t34} contains
the proofs of Theorems \ref{t:main-mb3} and \ref{t:main-mb4}. In Section \ref{s:discussion} we discuss some L\'evy processes  in $\R^d$ satisfying our assumptions.

We will use the following conventions in this paper.
 $c, c_0,
c_1, c_2, \cdots$ stand for constants
whose values are unimportant and which may change from one
appearance to another. All constants are positive finite numbers.
The labeling of the constants $c_0, c_1, c_2, \cdots$ starts anew in
the statement of each result. We will use ``$:=$"
to denote a definition, which is  read as ``is defined to be".
We denote $a \wedge b := \min \{ a, b\}$, $a \vee b := \max \{ a, b\}$.
Further, $f(t) \sim g(t)$, $t \to 0$ ($f(t) \sim g(t)$, $t \to
\infty$, respectively) means $ \lim_{t \to 0} f(t)/g(t) = 1$
($\lim_{t \to \infty} f(t)/g(t) = 1$, respectively).
Throughout the paper we will adopt the convention that
$X_{\zeta}=\partial$ and $u(\partial)=0$ for every function $u$.

\section{Oscillation reductions for inaccessible points}\label{s:inaccessible}

To handle the oscillation reductions at inaccessible points,
in this section we will assume, in addition to the corresponding assumptions in \cite{KSVp2},
that {\bf E1}$(z_0, R)$ ({\bf E2}$(z_0, R)$ respectively) holds
when we deal with finite boundary points (respectively infinity).

\subsection{Infinity}\label{ss:inaccessible-infty}

Throughout this subsection we will  assume that
$R_0=\infty$ and $D\subset \X$ is an unbounded open set.
We will deal with oscillation reduction at $\infty$ when $\infty$ is inaccessible from $D$ with respect $X$.
We further assume that there exists a point $z_0\in\X$
such that {\bf E2}$(z_0, R)$ and {\bf C2}$(z_0, R)$ hold, and that
$\widehat{X}$ satisfies \textbf{F2}$(z_0, R)$ for some $R>0$.
We will fix $z_0$ and $R$ and use the notation $B_r=B(z_0, r)$.
The next lemma is a direct consequence of assumption {\bf E2}$(z_0, R)$.

 \begin{lemma}\label{l:levy-density-I}
 For any $q\ge2$, $r\ge R$ and $\epsilon>0$,
  there exists $p=p(\epsilon,q,r)>16q$
  such that for every $z\in \overline{B}_{pr/8}^c$ and every $y\in \overline{B}_{qr}$, it holds that
 \begin{equation}\label{e:levy-density-I}
 (1+\epsilon)^{-1} < \frac{j(z,y)}{j(z, z_0)} < 1+\epsilon .
 \end{equation}
 \end{lemma}

In the remainder of this subsection, we assume that $r\ge R$,  and that $D$ is an open set such that $D\subset \overline{B}^c_r$. For $p>q>0$, let
$$
D^p=D\cap \overline{B}_p^c,\qquad D^{p,q}=D^q\setminus D^p.
$$

For $p>q>1$ and a nonnegative function $f$ on $\X$  define
 \begin{eqnarray}\label{f^pr}
 f^{pr,qr}(x)&=&
 \E_{x} \left[f(\wh X_{\wh \tau_{D^{pr}}}): \wh X_{\wh \tau_{D^{pr}}} \in D^{pr,qr}\right], \nonumber
\\
 \wt{f}^{pr,qr}(x)&=&
 \E_{x} \left[f(\wh X_{\wh \tau_{D^{pr}}}): \wh X_{\wh \tau_{D^{pr}}} \in (D\setminus D^{qr})\cup \overline{B}_r\right].
 \end{eqnarray}

\begin{lemma}\label{l:two-sided-estimate-I}
Suppose that $r\ge R$, $D\subset \overline{B}_r^c$ is an open set
and $f$ is a nonnegative function on $ \X$ which is  regular harmonic in $D$ with respect to $\wh X$
and vanishes on
$\overline{B}_r^c \cap ( \overline{D}^c \cup \wh D^{\mathrm{reg}})$.
 Let
$q\ge 2$, $\epsilon >0$, and choose $p=p(\epsilon, q,r)$ as in Lemma \ref{l:levy-density-I}.
Then for every $x\in D^{pr/8}$,
\begin{equation}\label{el:two-sided-estimate-I}
(1+\epsilon)^{-1} \wh P_{D^{pr/8}}(x,z_0) \int_{ \overline{B}_{qr}} f(y) m(dy)
\le \wt{f}^{pr/8,qr}(x) \le (1+\epsilon)
\wh P_{D^{pr/8}}(x,z_0) \int_{ \overline{B}_{qr}} f(y) m(dy).
\end{equation}
\end{lemma}
\pf Let $x\in D^{pr/8}$. Using Lemma \ref{l:levy-density-I} in the second inequality below, we get
\begin{eqnarray*}
\wt{f}^{pr/8,qr}(x) &= & \int_{D\setminus D^{qr}}  \wh P_{D^{pr/8}}(x,y)f(y)m(dy) +
\int_{\overline{B}_r} \wh P_{D^{pr/8}}(x,y)f(y)m(dy) \\
&=&\int_{D\setminus D^{qr}}\int_{D^{pr/8}} \wh G_{D^{pr/8}}(x,z) \wh j(z,y)m(dz) f(y)m(dy) \\
& & +\int_{\overline{B}_r}\int_{D^{pr/8}}\wh G_{D^{pr/8}}(x,z) \wh j(z,y)m(dz) f(y)m(dy) \\
&\le &(1+\epsilon)\int_{D\setminus D^{qr}}\int_{D^{pr/8}} \wh G_{D^{pr/8}}(x,z) \wh j(z, z_0)m(dz) f(y)m(dy) \\
& & +(1+\epsilon)\int_{\overline{B}_r}\int_{D^{pr/8}}\wh G_{D^{pr/8}}(x,z) \wh j(z, z_0)m(dz) f(y)m(dy) \\
&=&(1+\epsilon)\left(\int_{D\setminus D^{qr}} \wh P_{D^{pr/8}}(x,z_0)f(y)m(dy) +
\int_{\overline{B}_r}\wh P_{D^{pr/8}}(x,z_0)f(y)m(dy)\right) \\
&=&(1+\epsilon)\wh P_{D^{pr/8}}(x,z_0) \int_{ \overline{B}_{qr}} f(y) m(dy).
\end{eqnarray*}
This proves the right-hand side inequality. The left-hand side inequality can be proved in the same way. \qed

In the remainder of this subsection, we assume that
$r \ge R$, $D\subset \overline{B}_r^c$ an open set and $f_1$ and $f_2$ are nonnegative functions on $\X$ which are regular harmonic in $D$ with respect $\widehat{X}$ and vanish  on
$\overline{B}_r^c \cap ( \overline{D}^c \cup \wh D^{\mathrm{reg}})$.
Note that $f_i=f^{pr,qr}_i+ \wt{f}^{pr,qr}_i$.

\begin{lemma}\label{l:oscillation-assumption-1-I}
 Let $r\ge R$, $q>2$, $\epsilon>0$, and choose $p=p(\epsilon,q,r)$ as in
 Lemma \ref{l:levy-density-I}. If
 \begin{equation}\label{e:assumption-1-12-I}
\int_{ D^{3pr/8,qr}} f_i(y) m(dy)\le \epsilon \int_{ \overline{B}_{qr}} f_i(y) m(dy), \quad i=1,2,
 \end{equation}
 then, for all $x\in D^{pr}$.
\begin{equation}\label{e:estimate-of-quotient-I}
\frac{(1+\epsilon)^{-1}\int_{ \overline{B}_{qr}} f_1(y) m(dy)}{(C\epsilon +1+\epsilon)\int_{ \overline{B}_{qr}} f_2(y) m(dy)}\le \frac{f_1(x)}{f_2(x)} \le \frac{(C\epsilon +1+\epsilon)\int_{ \overline{B}_{qr}} f_1(y) m(dy)}{(1+\epsilon)^{-1}\int_{ \overline{B}_{qr}} f_2(y) m(dy)}.
\end{equation}
 \end{lemma}
 \pf Assume that  $x\in D^{pr}$.
Since  ${f_i}^{pr/8,qr}$ is regular harmonic in $D^{pr/8}$ with respect to $\wh{X}$ and vanishes
on
$\overline{B}_{pr/8}^c \cap ( \overline{D}^c \cup \wh D^{\mathrm{reg}})$,
using  {\bf F2}$(z_0, R)$ (with $a=3/2$),
 we have
 $$
f_i^{pr/8,qr}(x) \le C \wh{P}_{D^{pr/8}}(x,z_0)
  \int_{ B_{3pr/8}} f_i^{pr/8,qr}(y) m(dy).
 $$
 Since $f_i^{pr/8,qr}(y)\le f_i(y)$
  and $f_i^{pr/8,qr}(y)=0$ on $(D^{qr})^c$ except possibly at irregular points of $D$, by using that $m$ does not charge polar sets and applying \eqref{e:assumption-1-12-I} we have
 $$
 f_i^{pr/8,qr}(x) \le C \wh{P}_{D^{pr/8}}(x,z_0)
\int_{ D^{3pr/8,qr}} f_i(y) m(dy) \le C \epsilon  \wh{P}_{D^{pr/8}}(x,z_0)
   \int_{ \overline{B}_{qr}} f_i(y) m(dy).
 $$
 By  this and Lemma \ref{l:two-sided-estimate-I} we have
\begin{eqnarray*}
f_i(x)&=&f_i^{pr/8,qr}(x)+\wt{f}_i^{pr/8,qr}(x)\\
&\le & C\epsilon \wh{P}_{D^{pr/8}}(x,z_0)\int_{ \overline{B}_{qr}} f_i(y) m(dy)
+ (1+\epsilon) \wh{P}_{D^{pr/8}}(x,z_0)\int_{ \overline{B}_{qr}} f_i(y) m(dy)\\
&=&(C\epsilon +1+\epsilon)\wh{P}_{D^{pr/8}}(x,z_0)\int_{ \overline{B}_{qr}} f_i(y) m(dy)
\end{eqnarray*}
and
$$
f_i(x)\ge \wt{f}_i^{pr/8,qr}(x)\ge (1+\epsilon)^{-1}\wh{P}_{D^{pr/8}}(x,z_0)\int_{ \overline{B}_{qr}} f_i(y) m(dy).
$$
Therefore, \eqref{e:estimate-of-quotient-I} holds. \qed

Suppose that $\infty$ is inaccessible from $D$ with respect to $X$. Then there exists a point $x_0\in D$ such that
\begin{align}\label{e:pdRinfty}
\int_D G_D(x_0,y) m(dy) =\E_{x_0}\tau_D <\infty.
\end{align}
In the next result we fix this point $x_0$.

\begin{thm}\label{t:oscillation-reduction-yI}
Suppose that $\infty$ is inaccessible from  $D$ with respect to $X$.
Let $r > 2d(z_0,x_0) \vee R$.
For any two nonnegative functions $f_1$, $f_2$ on $\X$ which are regular harmonic in $D^r$ with respect to $\widehat{X}$ and vanish  on
$\overline{B}_r^c \cap ( \overline{D}^c \cup \wh D^{\mathrm{reg}})$
we have
\begin{equation}\label{e:inaccessible-limit-I}
\lim_{D\ni x\to \infty}\frac{f_1(x)}{f_2(x)}=
\frac{\int_{ \X} f_1(y) m(dy)}{\int_{ \X} f_2(y) m(dy)}.
\end{equation}
\end{thm}

\pf
First note that
$$
\int_{B_{3r}} G_D(x_0, z)m(dz) \ge \E_{x_0}[D \cap B_{3r}]>0.
$$
 By using  {\bf F2}$(z_0, R)$ we see that
$\int_{{B}_{8r}}f_i(y)m(dy)<\infty$.
The function $v\mapsto G_D(x_0,v)$ is regular harmonic in $D^r$ with respect to $\widehat{X}$ and vanishes on
$\overline{B}_{r}^c\setminus {D^{r}}$ (so vanishes on
$\overline{B}_r^c \cap ( \overline{D}^c \cup \wh D^{\mathrm{reg}})$).
By using  {\bf F2}$(z_0, R)$ for $\widehat{X}$, we have for $i=1, 2$,
\begin{align*}
&\int_{D^{8r}} f_i(y)m(dy) \le C
\int_{B_{3r}} f_i(z)m(dz)
\int_{D^{8r}} \wh{P}_{D^r} (y,z_0) m(dy)\\
&=C\int_{B_{3r}} G_D(x_0, z)m(dz) \int_{D^{8r}} \wh{P}_{D^r} (y,z_0) m(dy)
\frac{\int_{B_{3r}} f_i(z)m(dz)}{\int_{B_{3r}} G_D(x_0, z)m(dz)}
\\
& \le C^2
\int_{D^{8r}} G_D(x_0, y) m(dy)
 \frac{
\int_{B_{3r}} f_i(z)m(dz)}{\int_{B_{3r}} G_D(x_0, z)m(dz)} \\
& \le C^2
\int_{D} G_D(x_0, y) m(dy)
 \frac{
\int_{B_{3r}} f_i(z)m(dz)}{\E_{x_0}[D \cap B_{3r}]} < \infty.
\end{align*}
Hence $\int_{\X} f_i(y) m(dy)<\infty$, $i=1, 2$.
Let $q_0=2$ and $\epsilon >0$. For $j=0,1,\dots $, inductively define the
sequence $q_{j+1}=3p(\epsilon,q_j,r)/8 >6 q_{i}$ using Lemma \ref{l:levy-density-I}. Then
for $i=1, 2$,
$$
\sum_{j=0}^{\infty} \int_{ D^{q_{j+1}r, q_j r}} f_i(y) m(dy) =
\int_{ D^{q_0 r}} f_i(y) m(dy)<\infty .
$$
If $\int_{ D^{q_{j+1}r, q_j r}} f_i(y) m(dy)>\epsilon
\int_{ \overline{B}_{q_j r}} f_i(y) m(dy)$
for all $j\ge 0$, then
$$
\sum_{j=0}^{\infty}
\int_{ D^{q_{j+1}r, q_j r}} f_i(y) m(dy)
 \ge \epsilon
\sum_{j=0}^{\infty}
\int_{ \overline{B}_{q_j r}} f_i(y) m(dy)
\ge \epsilon \sum_{j=0}^{\infty}
\int_{ \overline{B}_{q_0 r}} f_i(y) m(dy)=\infty .
$$
Hence, there exists $k\ge 0$ such that
$\int_{ D^{q_{k+1}r, q_k r}} f_i(y) m(dy) \le \epsilon
\int_{ \overline{B}_{q_k r}} f_i(y) m(dy)$.
Moreover,
since $\lim_{j\to\infty}\int_{ D^{q_{j+1}r, q_j r}} f_i(y) m(dy)=0$, there exists $j_0\ge 0$
such that $\int_{ D^{q_{j+1}r, q_j r}} f_i(y) m(dy)\le
\int_{ D^{q_{k+1}r, q_k r}} f_i(y) m(dy)
$ for all $j\ge j_0$. Hence for all $j\ge j_0\vee k$ we have
$$
\int_{ D^{q_{j+1}r, q_j r}} f_i(y) m(dy)\le \int_{ D^{q_{k+1}r, q_k r}} f_i(y) m(dy) \le \epsilon
\int_{ \overline{B}_{q_k r}} f_i(y) m(dy) \le
\epsilon \int_{ \overline{B}_{q_j r}} f_i(y) m(dy).
$$
Therefore, there exists $j_0\in \N$ such that for all $j\ge j_0 \vee k$,
$$
\int_{ D^{q_{j+1}r, q_j r}} f_i(y) m(dy)\le
\epsilon \int_{ \overline{B}_{q_j r}} f_i(y) m(dy) \qquad i=1,2,
$$
and
$$
(1+\epsilon)^{-1}
\int_{ \X} f_i(y) m(dy)
 <
 \int_{ \overline{B}_{q_j r}} f_i(y) m(dy)
  <(1+\epsilon) \int_{ \X} f_i(y) m(dy),\qquad i=1,2.
$$
We see that the assumption of Lemma \ref{l:oscillation-assumption-1-I} are satisfied and conclude that \eqref{e:estimate-of-quotient-I} holds true: for $x\in D^{8q_{j+1}r/3}$,
$$
\frac{(1+\epsilon)^{-1} \int_{ \overline{B}_{q_j r}} f_1(y) m(dy)}{(C\epsilon +1+\epsilon) \int_{ \overline{B}_{q_j r}} f_2(y) m(dy)} \le
\frac{f_1(x)}{f_2(x)} \le \frac{(C\epsilon +1+\epsilon)\int_{ \overline{B}_{q_j r}} f_1(y) m(dy)}{(1+\epsilon)^{-1}
\int_{ \overline{B}_{q_j r}} f_2(y) m(dy)}.
$$
It follows that for $x\in D^{8q_{j+1}r/3}$,
$$
\frac{(1+\epsilon)^{-2}\int_{ \X} f_1(y) m(dy)}{(C\epsilon +1+\epsilon)(1+\epsilon)
\int_{ \X} f_2(y) m(dy)}\le \frac{f_1(x)}{f_2(x)} \le \frac{(C\epsilon +1+\epsilon)
(1+\epsilon)\int_{ \X} f_1(y) m(dy)}{(1+\epsilon)^{-2}\int_{ \X} f_2(y) m(dy)}.
$$
Since $\epsilon >0$ was arbitrary, we conclude that \eqref{e:inaccessible-limit-I} holds. \qed

\subsection{Finite boundary point}\label{ss:inaccessible-finite}

In this subsection, we deal with oscillation reduction at an inaccessible boundary point $z_0\in \X$ of an open set $D$.
Throughout the subsection, we assume that there exists $R\le R_0$
such that {\bf E1}$(z_0, R)$ and {\bf C1}$(z_0, R)$ hold,
and that $\widehat{X}$ satisfies \textbf{F1}$(z_0, R)$. We will fix this $z_0$.
Again, for simplicity, we use notation $B_r=B(z_0,r)$, $r>0$.

First, the next lemma is a direct consequence of assumption {\bf E1}$(z_0, R)$.

 \begin{lemma}\label{l:levy-density}
For any $q\in (0,1/2]$, $r\in (0,R]$ and $\epsilon>0$,
there exists $p=p(\epsilon,q,r)<q/16$ such that for every $z\in B_{8pr}$ and every $y\in B_{qr}^c$,
 \begin{equation}\label{e:levy-density}
 (1+\epsilon)^{-1} < \frac{j(z,y)}{j(z_0, y)} < 1+\epsilon .
 \end{equation}
 \end{lemma}

 Let $D\subset \X$ be an open set.
For $0<p<q$, let $D_p=D\cap B_p$ and $D_{p,q}=D_q\setminus D_p$.
For a function $f$ on $\X$, and $0<p<q$, let
 \begin{equation}\label{e:Lambda}
 \wh{\Lambda}_p(f):=\int_{\overline{B}_p^c} \wh j(z_0,y) f(y) m(dy),\qquad \wh{\Lambda}_{p,q}(f)=:\int_{D_{p,q}}\wh j(z_0,y) f(y) m(dy).
 \end{equation}
For $0<p<q<1$ and $r\in (0,R]$, define
 \begin{eqnarray}\label{f_pr}
 f_{pr,qr}(x)&=&
  \E_{x} \left[f(\wh X_{\wh \tau_{D_{pr}}}): \wh X_{\wh \tau_{D_{pr}}} \in D_{pr,qr}\right],\nonumber \\
 \wt{f}_{pr,qr}(x)&=&
 \E_{x} \left[f(\wh X_{\wh \tau_{D_{pr}}}): \wh X_{\wh \tau_{D_{pr}}} \in (D\setminus D_{qr})\cup B_r^c\right].
 \end{eqnarray}

\begin{lemma}\label{l:two-sided-estimate}
Let $q\in (0,1/2]$, $R\in (0,r]$, $\epsilon >0$, and choose $p=p(\epsilon, q,r)$
as in Lemma \ref{l:levy-density}.
Then for every $r\in (0,R]$, $D\subset B_r=B(z_0,r)$,  nonnegative function $f$ on $\X$ which is regular harmonic in $D$ with respect to $\widehat{X}$
and vanishes  on
$B_r \cap ( \overline{D}^c \cup\wh D^{\mathrm{reg}})$,
 and
every $x\in D_{8pr}$,
\begin{equation}\label{el:two-sided-estimate}
(1+\epsilon)^{-1}(\E_x \wh \tau_{D_{8pr}}) \wh{\Lambda}_{qr}(f)
\le \wt{f}_{8pr,qr}(x) \le (1+\epsilon)(\E_x  \wh \tau_{D_{8pr}}) \wh{\Lambda}_{qr}(f).
\end{equation}
\end{lemma}
\pf Let $x\in D_{8pr}$. Using Lemma \ref{l:levy-density} in the second inequality below, we get
\begin{eqnarray*}
\wt{f}_{8pr,qr}(x) &= & \int_{D\setminus D_{qr}}\wh P_{D_{8pr}}(x,y)f(y)m(dy) +\int_{B_R^c}\wh P_{D_{8pr}}(x,y)f(y)m(dy) \\
&=&\int_{D\setminus D_{qr}}\int_{D_{8pr}}\wh G_{D_{8pr}}(x,z) \wh j(z,y)m(dz) f(y)m(dy) \\
& & +\int_{B_R^c}\int_{D_{8pr}} \wh G_{D_{8pr}}(x,z)\wh j(z,y)m(dz) f(y)m(dy) \\
& \le & (1+\epsilon) (\E_x \wh \tau_{D_{8pr}}) \left(\int_{D\setminus D_{qr}} \wh j(z_0, y)f(y)m(dy) +\int_{B_r^c}\wh j(z_0, y)f(y)m(dy)\right)\\
&=& (1+\epsilon) (\E_x \wh \tau_{D_{8pr}})  \int_{B_{qr}^c}\wh j(z_0, y) f(y)m(dy)\\
&=&(1+\epsilon)( \E_x \wh \tau_{D_{8pr}}) \wh{\Lambda}_{qr}(f).
\end{eqnarray*}
This proves the right-hand side inequality. The left-hand side inequality can be
proved in the same way. \qed

In the remainder of this subsection, we assume $r\in (0,R]$, $D\subset B_r$ is an open set and $z_0\in \partial D$.
We also assume that $f_1$ and $f_2$ are nonnegative functions on $\X$ which are regular harmonic in $D$ with respect
to the process $\widehat{X}$, and vanish  on
$B_r \cap ( \overline{D}^c \cup \wh D^{\mathrm{reg}})$.
Note that $f_i=(f_i)_{pr,qr}+(\wt{f}_i)_{pr,qr}$.

 \begin{lemma}\label{l:oscillation-assumption-1}
 Let $R\in (0,1]$, $q<1/2$, $\epsilon>0$, and let $p=p(\epsilon,q,r)$ be as
 in Lemma \ref{l:levy-density}. If
 \begin{equation}\label{e:assumption-1-12}
 \wh{\Lambda}_{8pr/3,qr}(f_i)\le \epsilon \wh{\Lambda}_{qr}(f_i), \quad i=1,2,
 \end{equation}
 then for $x \in D_{pr}$
 \begin{equation}\label{e:estimate-of-quotient}
\frac{(1+\epsilon)^{-1}\wh{\Lambda}_{qr}(f_1)}{(C\epsilon +1+\epsilon)\wh{\Lambda}_{qr}(f_2)}\le \frac{f_1(x)}{f_2(x)} \le \frac{(C\epsilon +1+\epsilon)\wh{\Lambda}_{qr}(f_1)}{(1+\epsilon)^{-1}\wh{\Lambda}_{qr}(f_2)}.
\end{equation}
 \end{lemma}
 \pf Assume that $x\in D_{pr}$.
  Since  $(f_i)_{8pr,qr}$ is regular harmonic in $D_{8pr}$ with respect to $\wh X$ and vanish
on
$B_{8pr} \cap ( \overline{D}^c \cup \wh D^{\mathrm{reg}})$,
using  {\bf F1}$(z_0, R)$ (with $a=2/3$),
 we have
 $$
 (f_i)_{8pr,qr}(x) \le C (\E_x \tau_{D_{8pr}})
  \wh{\Lambda}_{8pr/3}( (f_i)_{8pr,qr}).
 $$
 Since $ (f_i)_{8pr,qr}(y)\le f_i(y)$
  and $ (f_i)_{8pr,qr}(y)=0$ on $D_{qr}^c$ except possibly at irregular points of $D$, applying \eqref{e:assumption-1-12} we have
 $$
 (f_i)_{8pr,qr}(x) \le C (\E_x \tau_{D_{8pr}})
  \wh{\Lambda}_{8pr/3,qr}(f_i) \le C \epsilon  (\E_x \tau_{D_{8pr}})
  \wh{\Lambda}_{qr}(f_i).
 $$
 By  this and Lemma \ref{l:two-sided-estimate},  we have that
\begin{eqnarray*}
f_i(x)&=&(f_i)_{8pr,qr}(x)+(\wt{f}_i)_{8pr,qr}(x)\\
&\le & C\epsilon (\E_x \tau_{D_{8pr}}) \wh{\Lambda}_{qr}(f_i)+ (1+\epsilon) (\E_x \tau_{D_{8pr}}) \wh{\Lambda}_{qr}(f_i)\\
&=&(C\epsilon +1+\epsilon) (\E_x \tau_{D_{8pr}}) \wh{\Lambda}_{qr}(f_i)
\end{eqnarray*}
and
$$
f_i(x)\ge (\wt{f}_i)_{8pr,qr}(x)\ge (1+\epsilon)^{-1}(\E_x \tau_{D_{8pr}}) \wh{\Lambda}_{qr}(f_i).
$$
Therefore,
\eqref{e:estimate-of-quotient} holds.
\qed

Assume that  $z_0$ is inaccessible  from  $D$ with respect to $X$.
Then there exist a point $x_0$ in $D$ such that
$$
P_D(x_0,z_0)=\int_D G_D(x_0,v)j(v, z_0) m(dv) <\infty .
$$
In the next result we fix this point $x_0$

\begin{thm}\label{t:reduction}
Suppose that $z_0$ is inaccessible from  $D$ with respect to $X$.
Let $r < 2d(z_0,x_0) \wedge R$.
For any two nonnegative functions $f_1$, $f_2$ on $\X$ which are regular harmonic in $D_r$ with respect to $\widehat{X}$ and vanish  on
$B_r \cap ( \overline{D}^c \cup\wh D^{\mathrm{reg}})$,
we have
\begin{equation}\label{e:inaccessible-limit}
\lim_{D\ni x\to z_0}\frac{f_1(x)}{f_2(x)}
=\frac{\int_{\X}\wh j(z_0, y)f_1(y)m(dy)}{\int_{\X}\wh j(z_0, y) f_2(y)m(dy)}.
\end{equation}
\end{thm}

\pf
 First note that
 $$
 \int_{\overline{B}_{r/3}^c} \wh j(z_0,z) G_D(x_0, z)m(dz) \ge \int_{D \cap \overline{B}_{r/3}^c} j(z, z_0) G_{D \cap \overline{B}_{r/3}^c} (x_0, z)m(dz)
 =P_{D \cap \overline{B}_{r/3}^c}(x_0,z_0) >0.
$$
Since $\widehat{X}$ satisfies {\bf F1}$(z_0, R)$, we have $\wh{\Lambda}_{r/8}(f_i)<\infty$.
The function $v\mapsto G_D(x_0,v)$ is regular harmonic in $D_r$ with respect to $\widehat{X}$ and vanishes on
$B_{r}\setminus D_{r}$ (so  vanishes on $B_r \cap ( \overline{D}^c \cup \wh D^{\mathrm{reg}})$).
By using  {\bf F1}$(z_0, R)$ for $\widehat{X}$ we have
\begin{align*}
&
\int_{B_{r/8}}\wh j(z_0, y)f_i(y) m(dy) \le C
\int_{\overline{B}_{r/3}^c}  \wh j(z_0,z) f_i(z)m(dz)
\int_{B_{r/8}}\wh j(z_0, y) \E_y [ \wh{\tau}_{D_r}] m(dy)
\\
&=C
\int_{\overline{B}_{r/3}^c} \wh j(z_0,z) G_D(x_0, z)m(dz)
\int_{B_{r/8}}\wh j(z_0, y) \E_y [ \wh{\tau}_{D_r}] m(dy)
\frac{\int_{\overline{B}_{r/3}^c}  \wh j(z_0,z) f_i(z)m(dz)}
{\int_{\overline{B}_{r/3}^c} \wh j(z_0,z) G_D(x_0, z)m(dz)}\\
&\le C^2 \int_{B_{r/8}}\wh j(z_0, y)G_D(x_0, y) m(dz)
\frac{\int_{\overline{B}_{r/3}^c}  \wh j(z_0,z) f_i(z)m(dz)}
{\int_{\overline{B}_{r/3}^c} \wh j(z_0,z) G_D(x_0, z)m(dz)}\\
&\le C^2 P_D(x_0,z_0)
\frac{ \wh{\Lambda}_{r/3}(f_i)   }
{P_{D \cap \overline{B}_{r/3}^c}(x_0,z_0) }<\infty.
\end{align*}
Therefore
\begin{eqnarray*}
\wh{\Lambda}(f_i):=\int_{\X}\wh j(z_0, y)f_i(y) m(dy)
=\int_{B_{r/8}}\wh j(z_0, y)f_i(y) m(dy)+\wh{\Lambda}_{r/8}(f_i)<\infty.
\end{eqnarray*}
Let $q_0=1/2$ and $\epsilon >0$. For $j=0,1,\dots $, inductively define the sequence
$q_{j+1}=p(\epsilon,q_j,r)$ as in Lemma \ref{l:levy-density}. Then
$$
\sum_{j=0}^{\infty} \wh{\Lambda}_{q_{j+1}r, q_j r}(f_i) =
\int_{D_{r/2}}\wh j(z_0, y)f_i(y) m(dy)\le \int_{\X}\wh j(z_0, y)f_i(y) m(dy)<\infty .
$$
If $\wh{\Lambda}_{q_{j+1}r, q_j r}(f_i)>\epsilon \wh{\Lambda}_{q_j r}(f_i)$
for all $j\ge 0$, then
$$
\sum_{j=0}^{\infty} \wh{\Lambda}_{q_{j+1}r, q_j r}(f_i) \ge \epsilon
\sum_{j=0}^{\infty} \wh{\Lambda}_{q_j r}(f_i)\ge \epsilon \sum_{j=0}^{\infty}
\wh{\Lambda}_{q_0 r}(f_i)=\infty .
$$
Hence, there exists an integer $k\ge 0$ such that $\wh{\Lambda}_{q_{k+1}r, q_k r}(f_i)
\le \epsilon \wh{\Lambda}_{q_k r}(f_i)$. Moreover, since $\lim_{j\to \infty}
\wh{\Lambda}_{q_{j+1}r, q_j r}(f_i)=0$, there exists $j_0\ge 0$ such that
$\wh{\Lambda}_{q_{j+1}r, q_j r}(f_i)\le \wh{\Lambda}_{q_{k+1}r, q_k r}(f_i)$
for all $j\ge j_0$. Hence for all $j\ge j_0\vee k$ we have
$$
\wh{\Lambda}_{q_{j+1}r, q_j r}(f_i)\le \wh{\Lambda}_{q_{k+1}r, q_k r}(f_i)
\le \epsilon \wh{\Lambda}_{q_k r}(f_i) \le \epsilon \wh{\Lambda}_{q_j r}(f_i).
$$
Therefore  for all $j\ge j_0$,
$$
\wh{\Lambda}_{q_{j+1}r, q_j r}(f_i) \le \epsilon \wh{\Lambda}_{q_j r}(f_i),\qquad i=1,2,
$$
and
$$
(1+\epsilon)^{-1}\wh{\Lambda}(f_i) < \wh{\Lambda}_{q_j r}(f_i) <(1+\epsilon) \wh{\Lambda}(f_i),\qquad i=1,2.
$$
Hence the assumption of Lemma \ref{l:oscillation-assumption-1} are satisfied and
consequently \eqref{e:estimate-of-quotient} holds: for $x\in D_{q_{j+1}r}$,
$$
\frac{(1+\epsilon)^{-1}\wh{\Lambda}_{q_j r}(f_1)}{(C\epsilon +1+\epsilon)\wh{\Lambda}_{q_j r}(f_2)}\le \frac{f_1(x)}{f_2(x)} \le \frac{(C\epsilon +1+\epsilon)\wh{\Lambda}_{q_j r}(f_1)}{(1+\epsilon)^{-1}\wh{\Lambda}_{q_j r}(f_2)}.
$$
It follows that $x\in D_{q_{j+1}r}$,
$$
\frac{(1+\epsilon)^{-2}\wh{\Lambda}(f_1)}{(C\epsilon +1+\epsilon)(1+\epsilon)\wh{\Lambda}(f_2)}\le \frac{f_1(x)}{f_2(x)} \le \frac{(C\epsilon +1+\epsilon)(1+\epsilon)\wh{\Lambda}(f_1)}{(1+\epsilon)^{-2}\wh{\Lambda}(f_2)}.
$$
Since $\epsilon >0$ was arbitrary, we conclude that \eqref{e:inaccessible-limit} holds. \qed

\section{Proof of Theorems \ref{t:main-mb3} and \ref{t:main-mb4}}\label{s:t34}

Let $D$ be a Greenian open subset of $\X$. Fix $x_0\in D$ and define
$$
M_D(x, y):=\frac{G_D(x, y)}{G_D(x_0, y)}, \qquad x, y\in D,~y\neq x_0.
$$

Combining \cite[Lemmas 3.2 and 3.4]{KSVp2}
and our Theorems \ref{t:oscillation-reduction-yI} and \ref{t:reduction}
we have the following.
\begin{thm}\label{t:1-10}
(a) Suppose that {\bf E1}$(z_0, R)$ holds
and that $\widehat{X}$ satisfies {\bf F1}$(z_0, R)$.
Then
\begin{align}\label{e:martin-kernel1}
M_D(x,z_0):=\lim_{D\ni v\to z_0}\frac{G_D(x,v)}{G_D(x_0,v)}
\end{align}
exists and is finite.
In particular, if $z_0$ is inaccessible from $D$ with respect to $X$, then
\begin{align}\label{e:martin-kernel2}
M_D(x,z_0)=\frac{\int_{\X}\wh j(z_0, y)G_D(x, y) m(dy)}{\int_{\X}\wh j(z_0, y)G_D(x_0, y) m(dy)}=  \frac{P_D(x,z_0)}{P_D(x_0,z_0)}\, .
\end{align}

\noindent
(b)
Suppose that {\bf E2}$(z_0, R)$ holds and that $\widehat{X}$ satisfies {\bf F2}$(z_0, R)$.
Then for every $x\in D$ the limit
\begin{align}\label{e:martin-kernel3}
M_D(x,\infty):=\lim_{D\ni v\to \infty}\frac{G_D(x,v)}{G_D(x_0,v)}
\end{align}
exists and is finite.
In particular, if $\infty$ is inaccessible from $D$ with respect to $X$, then
\begin{align}\label{e:martin-kernel4}
M_D(x,\infty)=\frac{\E_x \tau_D}{\E_{x_0}\tau_D}\, .
\end{align}
\end{thm}

Since both $X^D$ and $\widehat{X}^D$ are strongly Feller, the process $X^D$ satisfies Hypothesis (B) in \cite{KW}. See \cite[Section 4]{KSVp2} for details.
Therefore $D$  has
a Martin boundary $\partial_M D$ with respect to $X^D$ satisfying the following properties:
\begin{description}
\item{(M1)} $D\cup \partial_M D$ is
a compact metric space (with the metric denoted by $d_M$);
\item{(M2)} $D$ is open and dense in $D\cup \partial_M D$,  and its relative topology coincides with its original topology;
\item{(M3)}  $M_D(x ,\, \cdot\,)$ can be uniquely extended  to $\partial_M D$ in such a way that
\begin{description}
\item{(a)}
$ M_D(x, y) $ converges to $M_D(x, w)$ as $y\to w \in \partial_M D$ in the Martin topology;
\item{(b)} for each $ w \in D\cup \partial_M D$ the function $x \to M_D(x, w)$  is excessive with respect to $X^D$;
\item{(c)} the function $(x,w) \to M_D(x, w)$ is jointly continuous on
$D\times ((D\setminus\{x_0\})\cup \partial_M D)$ in the Martin topology and
\item{(d)} $M_D(\cdot,w_1)\not=M_D(\cdot, w_2)$ if $w_1 \not= w_2$ and $w_1, w_2 \in \partial_M D$.
\end{description}
\end{description}

\medskip

\noindent
\textbf{Proof of Theorem \ref{t:main-mb3}:}
(a) Using Theorem \ref{t:1-10}(a), by the same argument as in the proof of \cite[Theorem 1.1(a)]{KSVp2},
we have that $\partial_M^{z_0} D$ consists of a single point.

\noindent
(b) If $z_0$ is accessible from $D$ with respect to $X$, then by \cite[Theorem 1.1 (b)]{KSVp2} the Martin kernel $M_D(\cdot,z_0)$ is minimal harmonic for $X^D$.

Assume that $z_0$ is inaccessible from $D$ with respect to $X$.
Since $x\mapsto P_D(x,z_0)$ is \emph{not} harmonic with respect to $X^D$, we conclude
from \eqref{e:martin-kernel3}
that the Martin kernel $M_D(\cdot,z_0)$ is not harmonic, and in particular, that $z_0$ is \emph{not} a minimal Martin boundary point.
\qed

\noindent
\textbf{Proof of Corollary \ref{c:main-mb3}}: (a)
Let $\Xi:\partial D\to \partial^f_MD$ so that $\Xi(z)$ is the unique Martin boundary
point associated with $z\in \partial D$.
Since every finite Martin boundary point is associated with some $z\in \partial D$, we see that $\Xi$ is onto.
We show now that $\Xi$ is 1-1. If not, there are $z, z'\in \partial D$, $z\neq z'$, such that $\Xi(z)=\Xi(z')=w$.
Then $M_D(\cdot, z)=M_D(\cdot, w)= M_D(\cdot, z')$. It follows from the proof of \cite[Corollary 1.2(a)]{KSVp2}
that $z$ and $z'$ can not be both accessible. If one of them, say $z$, is accessible and the other, $z'$, is
inaccessible, then we can not have $M_D(\cdot, z)=M_D(\cdot, z')$ since $M_D(\cdot, z)$ is harmonic while
$M_D(\cdot, z')$ is not. Now let's assume that both $z$ and $z'$ are inaccessible.
Then $M_D(\cdot,z)=\frac{P_D(\cdot,z)}{P_D(x_0,z)}$ and
$M_D(\cdot,z')=\frac{P_D(\cdot,z')}{P_D(x_0,z')}$.
From $M_D(\cdot, z)=M_D(\cdot, z')$ we deduce that
$$
P_D(x,z)P_D(x_0,z')=P_D(x,z')P_D(x_0,z), \qquad \text{for all }x\in D.
$$
By treating $P_D(x_0,z')$ and $P_D(x_0,z)$ as constants, the above equality can be written as
$$
\int_D G_D(x,y)j(y, z)m(dy) = c \int_D G_D(x,y) j(y, z') m(dy), \qquad \text{for all }x\in D.
$$
By the uniqueness principle for potentials, this implies that the measures $j(y, z)m(dy)$ and $c j(y, z') m(dy)$ are equal. Hence $j(y, z)=c j(y, z')$ for $m$-a.e.~$y\in D$.
But this is impossible (for example, let $y\to z$; then $j(y, z)\to \infty$, while $cj(y, z')$ stays bounded
because of {\bf C1}$(z, R)$). We conclude that $z=z'$.

(b) The proof of this part is exactly the same as that of \cite[Corollary 1.2(b)]{KSVp2}.
\qed

\noindent
\textbf{Proof of Theorem \ref{t:main-mb4}:}

\noindent
(a)  Using Theorem \ref{t:1-10}(b), by the same argument as in the proof of \cite[Theorem 1.3(a)]{KSVp2}, we have that
 $\partial_M^{\infty} D$ is a single point which we will denote by $\infty$.

\noindent
(b) If $\infty$ is inaccessible from $D$ with respect to $X$,
then by \cite[Theorem 1.2 (b)]{KSVp2} the Martin kernel $M_D(\cdot, \infty)$ is minimal harmonic for $X^D$.

Assume that $\infty$ is inaccessible from $D$ with respect to $X$.
Since the function $x\mapsto \E_x \tau_D =\int_D G_D(x,y)m(dy)$
is \emph{not} harmonic with respect
to $X^D$, by \eqref{e:martin-kernel3} we conclude that the Martin kernel $M_D(\cdot, \infty)$ is not harmonic, and
in particular, $\infty$ is \emph{not} a minimal Martin boundary point. \qed

\noindent
\textbf{Proof of Corollary \ref{c:finite-not-infinite}}: (a) In the same way as in the proof of \cite[Corollary 1.4(a)]{KSVp2} it suffices to show that it cannot happen that $M_D(\cdot,\infty)=M_D(\cdot, z)$ for any $z\in \partial D$. If both $\infty$ and $z$ are accessible, this was shown in the proof of \cite[Corollary 1.4(a)]{KSVp2}. If one of the two points is accessible and the other inaccessible, then clearly the two Martin kernels are different since one is harmonic while the other is not. Assume that both $\infty$ and $z$ are inaccessible. Then $M_D(\cdot,\infty)=\frac{\E_{\cdot} \tau_D}{\E_{x_0}\tau_D}$ and $M_D(\cdot, z)=\frac{P_D(\cdot, z)}{P_D(x_0,z)}$. Therefore,
$$
P_D(x_0,z)\E_x \tau_D =P_D(x,z) \E_{x_0}\tau_D, \qquad \text{for all }x\in D.
$$
By treating $\E_{x_0}\tau_D$ and $P_D(x_0,z)$ as constants, the above equality can be written as
$$
\int_D G_D(x,y)m(dy) = c \int_D G_D(x,y) j(y, z) m(dy), \qquad \text{for all }x\in D.
$$
By the uniqueness principle for potentials, this implies that the measures $m(dy)$ and $c j(y, z) m(dy)$ are equal.  Hence $1=c j(y, z)$ for $m$-a.e.~$y\in D$ which clearly contradicts {\bf C1}$(z, R)$.

(b) The proof of this part is exactly the same as that of \cite[Corollary 1.4(b)]{KSVp2}.

\section{ Examples}\label{s:discussion}
In this section we discuss several classes of L\'evy processes  in $\R^d$ satisfying our assumptions.

\subsection{Subordinate Brownian motions}
In this subsection we discuss subordinate Brownian motions in $\R^d$ satisfying our assumptions.

 We will list conditions  on subordinate Brownian motions one by one under which our assumptions hold true.

Let  $W=(W_t, \P_x)$ be a Brownian motion in $\R^d$,
 $S=(S_t)$ an independent driftless subordinator with Laplace exponent $\phi$ and define the subordinate Brownian motion $Y=(Y_t, \P_x)$ by $Y_t=W_{S_t}$. Let $j_Y$ denote the L\'evy density of $Y$.

The Laplace exponent   $\phi$ is a  Bernstein function with $\phi(0+)=0$.  Since $\phi$ has no drift part,  $\phi$ can be written in the form $$
\phi(\lambda)=\int_0^{\infty}(1-e^{-\lambda t})\,
\mu(dt)\, .
$$
Here  $\mu$ is a $\sigma$-finite measure on
$(0,\infty)$ satisfying
$
\int_0^{\infty} (t\wedge 1)\, \mu(dt)< \infty.
$
 $\mu$ is called the L\'evy measure
of the subordinator $S$.
$\phi$  is called a complete Bernstein function
if the L\'evy measure $\mu$ of $S_t$ has a completely monotone density
$\mu(t)$, i.e., $(-1)^n D^n \mu\ge 0$ for every non-negative integer
$n$.
We will assume that $\phi$ is a complete Bernstein function.

When $\phi$ is unbounded and $Y$ is transient,
the mean occupation time measure of $Y$ admits a density $G(x,y)=g(|x-y|)$
which is called the Green function of $Y$,
and is given by the formula
\begin{equation}\label{e:green-function}
g(r):=\int_0^{\infty}(4\pi t)^{-d/2} e^{-r^2/(4t)}u(t)\, dt\, .
\end{equation}
Here $u$ is the potential density of the subordinator $S$.

\bigskip
We first discuss conditions that ensure {\bf E1}$(z_0, R)$.

By \cite[Lemma A.1]{KM}, for all $t>0$, we have
\begin{equation}\label{e:muexp}
\mu(t)\le (1-2e^{-1})^{-1}t^{-2}\phi'(t^{-1}) \le  (1-2e^{-1})^{-1}t^{-1}\phi(t^{-1}).
\end{equation}
 Thus
\begin{equation}\label{e:expub4mu}
\mu(t)\le (1-2e^{-1})^{-1}\phi'(M^{-1})t^{-2}, \qquad t \in (0,  M].
\end{equation}

In \cite{KSV12b}, we have shown that there exists $c\in (0, 1)$ such that
\begin{equation}\label{e:bofmuatinfty}
\mu(t+1)\ge c\mu(t), \qquad t\ge 1.
\end{equation}
As a consequence of this, one can easily show that there exist $c_1, c_2>0$ such that
\begin{equation}\label{e:explb4mu}
\mu(t)\ge c_1e^{-c_2 t}, \qquad t\ge 1.
\end{equation}
In fact, it follows from \eqref{e:bofmuatinfty} that for any $n\ge 1$,
$
\mu(n+1)\ge c^n\mu(1).
$
Thus, for any $t\ge 1$,
\begin{eqnarray*}
\mu(t)&\ge&\mu([t]+1)\ge c^{[t]}\mu(1)=\mu(1)e^{[t]\log c}\\
&=&\mu(1)e^{([t]-t)\log c}e^{t\log c}\ge c^{-1}\mu(1)e^{t\log c}.
\end{eqnarray*}

The following is a refinement of \eqref{e:bofmuatinfty} and \cite[Lemma 3.1]{KL}.

\begin{lemma}\label{l:1}
Suppose that  the Laplace exponent
$\phi$ of $S$ is a complete Bernstein function.
Then, for any $t_0>0$,
$$
\lim_{\delta\to 0} \sup_{t>t_0}\frac{\mu(t)}{\mu(t+\delta)}=1\, .
$$
\end{lemma}

\pf This is proof is similar to the proof of \cite[Lemma 3.1]{KL}, which in turn
is a refinement of the proof of \cite[Lemma 13.2.1]{KSV12b}. Let $\eta>0$ be given. Since $\mu$ is
a complete monotone function, there exists a measure $m$ on $[0, \infty)$ such that
$$
\mu(t)=\int_{[0, \infty)}e^{-tx}m(dx), \qquad t>0.
$$
Choose $r=r(\eta, t_0)>0$ such that
$$
\eta\int_{[0, r]}e^{-t_0x}m(dx)\ge \int_{(r, \infty)}e^{-t_0x}m(dx).
$$
Then for any $t>t_0$, we have
\begin{eqnarray*}
\lefteqn{\eta\int_{[0, r]}e^{-tx}m(dx)=\eta\int_{[0, r]}e^{-(t-t_0)x}e^{-t_0x}m(dx)\ge \eta e^{-(t-t_0)r}\int_{[0, r]}e^{-t_0x}m(dx)}\\
&\ge& e^{-(t-t_0)r}\int_{(r, \infty)}e^{-t_0x}m(dx)=\int_{(r, \infty)}e^{-(t-t_0)r}e^{-t_0x}m(dx)\ge  \int_{(r, \infty)}e^{-tx}m(dx).
\end{eqnarray*}
Thus for any $t>t_0$ and $\delta>0$,
\begin{eqnarray*}
\mu(t+\delta)&\ge&\int_{[0, r]}e^{-(t+\delta)x}m(dx)\ge e^{-r\delta}\int_{[0, r]}e^{-tx}m(dx)\\
&=&e^{-r\delta}(1+\eta)^{-1}\left(\int_{[0, r]}e^{-tx}m(dx)+\eta\int_{[0, r]}e^{-tx}m(dx) \right)\\
&\ge&e^{-r\delta}(1+\eta)^{-1}\left(\int_{[0, r]}e^{-tx}m(dx)+\int_{(r, \infty)}e^{-tx}m(dx) \right)\\
&=&e^{-r\delta}(1+\eta)^{-1}\int_{[0, \infty)}e^{-tx}m(dx)=e^{-r\delta}(1+\eta)^{-1}\mu(t).
\end{eqnarray*}
Therefore
$$
\limsup_{\delta\to 0}\sup_{t>t_0}\frac{\mu(t)}{\mu(t+\delta)}\le 1+\eta.
$$
Since $\eta$ is arbitrary and $\mu$ is decreasing, the assertion of the lemma is valid.
\qed

 The L\'evy measure of $Y$ has a density with respect to
the Lebesgue measure given by $j_Y(x)=j(|x|)$ with
$$
j(r)=\int^\infty_0 g(t, r)\mu(t)dt, \qquad r\neq 0,
$$
where
$$
g(t, r)=(4\pi t)^{-d/2}\exp(-\frac{r^2}{4t}).
$$

As a consequence of \eqref{e:bofmuatinfty}, one can easily get that there exists
$c\in (0, 1)$ such that
\begin{equation}\label{e:bofjatinfty}
j(r+1)\ge cj(r), \qquad r\ge 1.
\end{equation}
Using this, we can show that there exist $c_1, c_2>0$ such that
\begin{equation}\label{e:explb4j}
j(r)\ge c_1e^{-c_2 r}, \qquad r\ge 1.
\end{equation}

\begin{lemma}\label{l:2}
Suppose that  the Laplace exponent
$\phi$ of $S$ is a complete Bernstein function.
For any $r_0\in (0, 1)$,
$$
\lim_{\eta\to 0} \sup_{r>r_0}\frac{\int^\eta_0g(t, r)\mu(t)dt}{j(r)}=0\, .
$$
\end{lemma}

\pf For any $\eta\in (0, 1)$ and $r\in (r_0, 2]$, we have
$$
\frac{\int^\eta_0g(t, r)\mu(t)dt}{j(r)}\le \frac{\int^\eta_0g(t, r_0)\mu(t)dt}{j(2)}.
$$
Thus
$$
\lim_{\eta\to 0}\sup_{r\in (r_0, 2]}\frac{\int^\eta_0g(t, r)\mu(t)dt}{j(r)}
=0.
$$
Thus we only need to show that
$$
\lim_{\eta\to 0}\sup_{r>2}\frac{\int^\eta_0g(t, r)\mu(t)dt}{j(r)}
=0.
$$

It follows from \eqref{e:expub4mu} that
\begin{eqnarray*}
\lefteqn{\int^\eta_0(4\pi t)^{-d/2}\exp(-\frac{r^2}{4t})\mu(t)dt
\le c_1\int^\eta_0t^{-(\frac{d}2+2)}\exp(-\frac{r^2}{4t})dt\le c_3\int^\eta_0\exp(-\frac{r^2}{8t})dt}\\
&=& c_3 \int^\infty_{r^2/(8\eta)}e^{-s}\frac{r^2}{8s^2}ds\le c_4r^2\int^\infty_{r^2/(8\eta)}e^{-s/2}ds=c_5r^2\exp(-\frac{r^2}{16\eta}).
\end{eqnarray*}
Now combining this with \eqref{e:explb4j} we immediately arrive at the
desired conclusion.
\qed

\begin{lemma}\label{l:3}
Suppose that  the Laplace exponent
$\phi$ of $S$ is a complete Bernstein function.
For any $r_0\in (0, 1)$,
\begin{align}\label{e:BL1}
\lim_{\delta\to 0} \sup_{r>r_0}\frac{j(r)}{j(r+\delta)}=1\,  .
\end{align}
\end{lemma}

\pf For any $\epsilon\in (0, 1)$, choose $\eta\in (0, 1)$ such that
$$
\sup_{r>r_0}\frac{\int^\eta_0 g(t, r)\mu(t)dt}{j(r)}\le \epsilon.
$$
Then for any $r>r_0$,
$
\int^\infty_\eta g(t, r)\mu(t)dt\ge (1-\epsilon)j(r).
$
Fix this $\eta$. It follows from Lemma \ref{l:1} that there exists $\delta_0
\in (0, \eta/2)$ such that
$$
\frac{\mu(t)}{\mu(t+\delta)}\le 1+\epsilon, \qquad t\ge \eta, \delta\in (0, \delta_0].
$$
For $t>\eta$, $0\le (r+\delta-t)^2=(r+\delta)^2-2tr+t(t-\delta)-\delta t$ and so
$t(t-\delta)\ge 2tr+\delta t-(r+\delta)^2$. Thus
\begin{eqnarray*}
\frac{(r+\delta)^2}{4t}-\frac{r^2}{4(t-\delta)}=
\frac{(r+\delta)^2(t-\delta)-r^2t}{4t(t-\delta)}
=\frac{\delta(2tr+\delta t-(r+\delta)^2)}{4t(t-\delta)}\le \frac\delta4.
\end{eqnarray*}
Consequently, for $r>r_0$ and $\delta\in (0, \delta_0)$,
\begin{eqnarray*}
j(r+\delta)&\ge & \int^\infty_\eta(4\pi t)^{-d/2}\exp(-\frac{(r+\delta)^2}{4t})\mu(t)dt\\
&\ge &e^{-\delta/4}\int^\infty_\eta (4\pi t)^{-d/2}\exp(-\frac{r^2}{4(t-\delta)})\mu(t)dt\\
&\ge &e^{-\delta/4}\int^\infty_{\eta-\delta} (4\pi (t+\delta))^{-d/2}\exp(-\frac{r^2}{4t})\mu(t+\delta)dt\\
&\ge &e^{-\delta/4}\left(\frac{\eta}{\eta+\delta}\right)^{d/2}(1+\epsilon)^{-1}\int^\infty_\eta g(t, r)\mu(t)dt\\
&\ge&e^{-\delta/4}\left(\frac{\eta}{\eta+\delta}\right)^{d/2}(1+\epsilon)^{-1}(1-\epsilon)j(r).
\end{eqnarray*}
Now choose $\delta^*\in (0, \delta_0)$ such that
$$
e^{-\delta/4}\left(\frac{\eta}{\eta+\delta}\right)^{d/2}\ge (1+\epsilon)^{-1},
\qquad \delta \in (0, \delta^*].
$$
Then for all $r>r_0$ and $\delta \in (0, \delta^*]$,
$$
j(r+\delta)\ge (1+\epsilon)^{-2}(1-\epsilon)j(r),
$$
which is equivalent to
$$
\frac{j(r)}{j(r+\delta)}\le \frac{(1+\epsilon)^2}{(1-\epsilon)},
$$
which implies \eqref{e:BL1}.
\qed

 \begin{lemma}\label{l:levy-densityn2}
If  the Laplace exponent
$\phi$ of $S$ is a complete Bernstein function,
then {\bf E1}$(z_0, R)$ holds for $Y$.
\end{lemma}
 \pf
Fix  $r_0, \eps>0$ and use the notation $B_r=B(0, r)$.
By Lemma \eqref{l:3}
there exists $\eta=\eta(\epsilon,r_0)>0$ such that for all $\eta\le \eta(\epsilon,r_0)$,
 $$
  \sup_{r>r_0} \frac{j(r)}{j(r+\eta)} <1+\epsilon .
 $$
 Let
 $
 \delta:=\frac{2\eta}{r_0}\wedge 1.
 $
  For $y\in B_{\delta r_0/2}$ and $z\in B_{2r_0}^c$ we have
 \begin{eqnarray*}
&&r_0 <\frac{|z|}{2}= |z|-\frac{|z|}{2} \le |z|-|y|  \le  |z-y|\le |z|+|y| \le |z|+\frac{\delta r_0}{2}  \le |z|+\eta,\\
&& r_0<|z| \le  |z-y|+|y|\le |z-y|+\eta .
 \end{eqnarray*}
 Hence,
 $$
  \frac{j(|z-y|)}{j(|z|)}\le \frac{j(|z-y|)}{j(|z-y|+\eta)}\le \sup_{r>r_0}
 \frac{j(r)}{j(r+\eta)}< 1+\epsilon
 $$
 and
 $$
  \frac{j(|z|)}{j(|z-y|)}\le \frac{j(|z|)}{j(|z|+\eta)} \le \sup_{r>r_0}
 \frac{j(r)}{j(r+\eta)}< 1+\epsilon.
 $$
This finishes the proof of the lemma. \qed

We now briefly discuss \eqref{e:G111},  {\bf C1}$(z_0, R)$, \eqref{e:nGassup1},   \textbf{F1}$(z_0, R)$,  and {\bf G}.
First note that,  if $Y$ is transient  then \eqref{e:G111} holds (see \cite[Lemma 2.10]{KSV14b}).
For the remainder of this section, we will always assume that
$\phi$ is a complete Bernstein function and  the L\'evy density $\mu$ of $\phi$ is infinite, i.e.
$\mu(0,\infty)=\infty$.
We consider the following further assumptions  on  $\phi$:

\medskip

\noindent
{\bf H}:
{\it  there exist constants $\sigma>0$, $\lambda_0 > 0$ and
$\delta \in (0, 1]$
 such that
\begin{equation}\label{e:sigma}
  \frac{\phi'(\lambda t)}{\phi'(\lambda)}\leq\sigma\, t^{-\delta}\ \text{ for all }\ t\geq 1\ \text{ and }\
 \lambda \ge \lambda_0\, .
\end{equation}
When $d \le 2$,  we assume that $d+2\delta-2>0$ where $\delta$ is the constant
in \eqref{e:sigma}, and  there are $\sigma'>0$ and
\begin{equation}\label{e:new22}
\delta'  \in  \left(1-\tfrac{d}{2}, (1+\tfrac{d}{2})
\wedge (2\delta+\tfrac{d-2}{2})\right)
\end{equation}
 such that
\begin{equation}\label{e:new23}
  \frac{\phi'(\lambda x)}{\phi'(\lambda)}\geq
\sigma'\,x^{-\delta'}\ \text{ for all
}\ x\geq 1\ \text{ and }\ \lambda\geq\lambda_0\,.
\end{equation}
}
Assumption {\bf H}  was introduced and used in \cite{KM} and \cite{KM2}. It is easy to check that
if $\phi$ is a complete Bernstein function satisfying
satisfying a weak lower scaling condition at  infinity
\begin{equation}\label{e:new2322}
a_1 \lambda^{\delta_1}\phi(t)\le \phi(\lambda t)\le a_2 \lambda^{\delta_2}\phi(t)\, ,\qquad \lambda \ge 1, t\ge 1\, ,
\end{equation}
for some $a_1, a_2>0$ and $\delta_1, \delta_2\in (0,1)$, then {\bf H}  is automatically satisfied.
One of the reasons for adopting the more general setup above is to cover the case
of geometric stable and iterated geometric stable subordinators.
Suppose that $\alpha\in (0, 2)$ for $d \ge 2$ and that $\alpha\in (0, 2]$ for $d \ge 3$.
A geometric $(\alpha/2)$-stable subordinator is a subordinator
with Laplace exponent $\phi(\lambda)=\log(1+\lambda^{\alpha/2})$.
Let $\phi_1(\lambda):=\log(1+\lambda^{\alpha/2})$, and for $n\ge 2$,
$\phi_n(\lambda):=\phi_1(\phi_{n-1}(\lambda))$. A subordinator with
Laplace exponent $\phi_n$ is called an iterated geometric subordinator.
It is easy to check that the functions $\phi$ and $\phi_n$ satisfy
{\bf H}
but they do not satisfy \eqref{e:new2322}.

It follows from \cite[Lemma 5.4]{{KM2}} and \cite[Section 4.2]{KSVp2} that if $Y$ is transient and {\bf H} is true, then there exists $R>0$ such that the assumptions  {\bf G}, {\bf C1}$(z_0, R)$,  \eqref{e:nGassup1} and \textbf{F1}$(z_0, R)$  hold for all $z_0 \in \R^d$.
Thus using these facts and Lemma \ref{l:levy-densityn2}, we have the following as a special case of Theorems \ref{t:main-mb0}(a) and  \ref{t:main-mb3}(b).

 \begin{corollary}\label{c:lK1}
 Suppose that
$Y=(Y_t, \P_x:\, t\ge 0, x \in \R^d)$
 is a transient subordinate Brownian motion
whose
characteristic exponent is given by $\Phi(\theta)=\phi(|\theta|^2)$,
$\theta\in \R^d$.
Suppose $\phi$ is a complete Bernstein function
with  the infinite L\' evy measure $\mu$
and assume that {\bf H} holds.
Let $r  \le 1$ and let $f_1$ and $f_2$ be nonnegative functions on $\R^d$ which are regular harmonic in $D\cap {B}(z_0, r)$ with respect to the process $Y$, and
vanish  on
$B(z_0, r) \cap ( \overline{D}^c \cup D^{\mathrm{reg}})$.
 Then the limit
$$
\lim_{D\ni x\to z_0}\frac{f_1(x)}{f_2(x)}
$$
exists and is finite. Moreover, the Martin boundary point associated with $z\in \partial D$ is minimal
if and only if $z$ is accessible from $D$.
 \end{corollary}

\subsection{Unimodal L\'evy process}

Let $Y$ be an  isotropic unimodal L\'evy process whose characteristic exponent is
$\Psi_0(|\xi|)$, that is,
\begin{equation}\label{e:fuku1.3}
\Psi_0(|\xi|)= \int_{\R^d}(1-\cos(\xi\cdot y))j_0(|y|)dy
\end{equation}
where the function $x \mapsto j_0(|x|)$ is the L\'evy density of $Y$.
If $Y$ is transient, let $x \mapsto g_0(|x|)$ denote the Green function of $Y$.

Let $0<\alpha < 2$. Suppose that $\Psi_0(\lambda)\sim \lambda^{\alpha}\ell(\lambda)$,
$\lambda \to 0$, and $\ell$ is a slowly varying function at $0$.
Then by \cite[Theorems 5 and 6]{CGT} we have  the following asymptotics of $j_0$ and $g_0$.
 \begin{lemma}\label{l:j-g-near-infty} Suppose that $\Psi_0(\lambda)\sim \lambda^{\alpha}\ell(\lambda)$,
$\lambda \to 0$, and $\ell$ is a slowly varying function at $0$ ,
\begin{itemize}
    \item[(a)] It holds that
    \begin{equation}\label{e:j-near-infty}
    j_0(r)\sim r^{-d}
    \Psi_0(r^{-1})\, ,\qquad r\to \infty.
    \end{equation}
    \item[(b)] If $d  \ge 3$, then $Y$ is transient and
    \begin{equation}\label{e:g-near-infty}
    g_0(r)
    \sim  r^{-d}
    \Psi_0(r^{-1})^{-1}\, ,\qquad r\to \infty.
    \end{equation}
   \end{itemize}
\end{lemma}
We further assume that the L\'{e}vy measure of $X$ is infinite. Then by \cite[Lemma 2.5]{KR} the density function $x \to p_t(|x|)$ of $X$ is continuous and, by the strong Markov property, so is the density function of $X^D$.
 Using the upper bound of $p_t(|x|)$ in \cite[Theorem 2.2]{KKK} (which works for all $t>0$) and the monotonicity of  $r \to p_t(r)$, we see that the Green function of $X^D$ is continuous  for all open set $D$.
From this and \eqref{e:g-near-infty}, we have that if $d \ge 3$, the L\'{e}vy measure is infinite, and  $\Psi_0(\lambda)\sim \lambda^{\alpha}\ell(\lambda)$, then  ${\bf G}$ and \eqref{e:nGassup1-infty} hold (see \cite[Proposition 6.2]{KSVp1}).
It is proved in \cite{KSVp1} under some assumptions
much weaker than the above that {\bf F2}$(z_0, R)$ holds for all $z_0 \in \R^d$.
From \eqref{e:j-near-infty} we have
that {\bf E2}$(z_0, R)$ and {\bf C2}$(z_0, R)$ hold for all $z_0 \in \R^d$.
Using the above facts, we have the following as a special case of Theorems \ref{t:main-mb0}(b) and  \ref{t:main-mb4}(b).
 \begin{corollary}\label{c:lK2}
 Suppose that $d \ge 3$ and that
$Y=(Y_t, \P_x:\, t\ge 0, x \in \R^d)$
 is an isotropic unimodal L\'evy process
whose
characteristic exponent is given by $\Psi_0(|\xi|)$.
Suppose that $0<\alpha < 2$, that the L\'{e}vy measure of $X$ is infinite, and that
$\Psi_0(\lambda)\sim \lambda^{\alpha}\ell(\lambda)$,
$\lambda \to 0$, and $\ell$ is a slowly varying function at $0$.
Let $r > 1$, $D$ be an unbounded open set and let $f_1$ and $f_2$ be nonnegative functions on $\R^d$ which are regular harmonic in $D\cap \overline{B}(z_0, r)^c
$ with respect to the process $Y$, and vanish  on
$\overline{B}(z_0, r)^c\cap ( \overline{D}^c \cup D^{\mathrm{reg}})$.
Then the limit
$$
\lim_{D\ni x\to \infty}\frac{f_1(x)}{f_2(x)}
$$
exists and is finite. Moreover,
the Martin boundary point
associated with $\infty$ is minimal if and only if $\infty$ is accessible from $D$.
 \end{corollary}

\begin{remark}\label{r:B0}
{\rm Using  \cite[Lemma 3.3]{KSV8} instead of \cite[Theorem 6]{CGT}, one can see that
Corollary \ref{c:lK2} holds for $d > 2\alpha$
when $Y$
 is a subordinate Brownian motion
whose
Laplace exponent $ \phi$ is a complete Bernstein function
and that  $\phi(\lambda)\sim \lambda^{\alpha/2}\ell(\lambda)$ where $0<\alpha < 2$ and $\ell$ is a slowly varying function at $0$,
}
\end{remark}

\begin{remark}\label{r:B1}
{\rm
If $Y$ is a L\'evy process satisfying {\bf E1}$(z_0, R)$, ({\bf E2}$(z_0, R)$, respectively)
then the L\'evy process $Z$ with Levy density $j_Z(x):=k(x/|x|)j_Y(x)$ also satisfies
{\bf E1}$(z_0, R)$ ({\bf E2}$(z_0, R)$, respectively) when $k$ is a continuous
function on the unit sphere and bounded between two positive constant.
In fact, since
$$
\left|\frac{z-y}{|z-y|}-\frac{z}{|z|}\right|
\le\left |\frac{z-y}{|z-y|} -\frac{z-y}{|z|} \right|+\left|\frac{z-y}{|z|} -\frac{z}{|z|}\right|
\le \frac{||z|-|z-y||}{|z|}+\frac{|y|}{|z|}
\le \frac{2|y|}{|z|},
$$
we have
$$
\left|\frac{z-y}{|z-y|}-\frac{z}{|z|}\right| \le \frac{2|y|}{r} \quad \text{for all }|z|>r \quad
\text{and} \quad
\left|\frac{z-y}{|z-y|}-\frac{z}{|z|}\right| \le \frac{2r}{|z|} \quad \text{for all }|y|<r.
$$
Moreover, since $k$ is bounded below by a positive constant,
$$
\left|\frac{k(z/|z|)}{k((z-y)/|z-y|)}-1 \right| \le c |{k(z/|z|)}-{k((z-y)/|z-y|)}|.
$$
Thus by uniform continuity of $k$ on the unit sphere, we see that for all $r>0$
$$
\lim_{|y|  \to 0}\sup_{z: |z|>r}\frac{ k(z/|z|) }{k((z-y)/|z-y|)}
=\lim_{|y|  \to 0}\sup_{z: |z|>r}\frac{ k(z/|z|) }{k((z-y)/|z-y|)}
=1,
$$
and
$$\lim_{|z| \to \infty} \sup_{y: |y| < r}\frac{ k(z/|z|) }{k((z-y)/|z-y|)}
=\lim_{|z| \to \infty} \inf_{y: |y| < r}
\frac{k(z/|z|) }{k((z-y)/|z-y|)}
=1.
$$
When $Y$ is a symmetric stable process, this includes not necessarily
symmetric strictly stable processes with Levy density $c k(x/|x|)|x|^{-d-\alpha}$ where $k$ is a continuous function on the unit sphere bounded between two
positive constant.
}
\end{remark}

\bigskip
\noindent
{\bf Acknowledgements:} Part of the research for this paper was done
during the visit of Renming Song and Zoran Vondra\v{c}ek to Seoul National University
from May 24 to June 8 of 2015.
They thank the Department of
Mathematical Sciences of Seoul National University for the hospitality.

\end{doublespace}

\bigskip

\vspace{.1in}
\begin{singlespace}


\end{singlespace}

\
\vskip 0.1truein

\parindent=0em

{\bf Panki Kim}

Department of Mathematical Sciences and Research Institute of Mathematics,

Seoul National University, Building 27, 1 Gwanak-ro, Gwanak-gu Seoul 08826, Republic of Korea

E-mail: \texttt{pkim@snu.ac.kr}

\bigskip

{\bf Renming Song}

Department of Mathematics, University of Illinois, Urbana, IL 61801,
USA

E-mail: \texttt{rsong@math.uiuc.edu}

\bigskip

{\bf Zoran Vondra\v{c}ek}

Department of Mathematics, University of Zagreb, Zagreb, Croatia, and \\
Department of Mathematics, University of Illinois, Urbana, IL 61801,
USA

Email: \texttt{vondra@math.hr}
\end{document}